\documentclass[letterpaper, 10 pt, conference]{ieeeconf}  % Comment this line out if you need a4paper

\IEEEoverridecommandlockouts                              % This command is only needed if 
\title{\LARGE \bf Graphical Newton}

\author{Akshay Srinivasan$^{1}$ and Emanuel Todorov$^{2}$
\thanks{$^{1}$ This work was done when Akshay Srinivasan was at the University of Washington {\tt akshaysrinivasan@gmail.com}}
\thanks{$^{2}$ Emanuel Todorov is with the Depts. of Computer Science \& Engineering and Applied Mathematics, University of Washington, Seattle, WA 98195, USA {\tt todorov@cs.washington.edu}}}
\usepackage{graphicx} % more modern
\usepackage{subfigure} 

\usepackage{algorithm}
\usepackage{algorithmic}

\usepackage{url}
\usepackage{enumerate}
\usepackage{amsmath}
\usepackage{amssymb}
\usepackage{amsfonts}
\usepackage{color}

\newcommand{\figref}[1]{(Figure~\ref{#1})}
\newcommand{\eqnref}[1]{(\ref{#1})}

\newcommand{\diff}[1]{{\mathtt{d} {#1}}}
\newcommand{\norm}[1]{\vert \vert {#1} \vert \vert}

\newcommand{\qed}{\(\hfill \square\)}

\newcommand{\p}{\partial}

\DeclareMathAlphabet{\mathpzc}{OT1}{pzc}{m}{it}
\newcommand{\obj}[1]{{\mathpzc{#1}}}
\newcommand{\pr}[1]{\mathtt{S}_{#1}}

\DeclareMathOperator{\supp}{supp}
\DeclareMathOperator{\rank}{rank}

\DeclareMathOperator{\e}{e}

\DeclareMathOperator{\tw}{tw}
\DeclareMathOperator{\pa}{\delta^{+}} %\delta^{-}} %pa
\DeclareMathOperator{\ch}{\delta^{-}}  %\delta^{+}} %ch
\DeclareMathOperator{\nh}{\delta}  %\delta^{+}} %ch

\newtheorem{lem}{Lemma}
\newtheorem{theo}{Theorem}
\newtheorem{defn}{Definition}
\usepackage[latin1]{inputenc}
\usepackage{tikz}
\usetikzlibrary{trees, arrows, matrix, positioning}

\begin{document}

\maketitle
\thispagestyle{empty}
\pagestyle{empty}

\begin{abstract}
  Computing the Newton step for a generic function $\obj{f}: \mathbb{R}^N \rightarrow \mathbb{R}$ takes $O(N^{3})$ flops. In this paper, we explore avenues for reducing this bound, when the computational structure of \(f\) is known beforehand. It is shown that the Newton step can be computed in time, linear in the size of the computational-graph, and cubic in its tree-width.
\end{abstract}

\section{Introduction}
Newton's method is essential to many areas of at the core of second-order methods in nonlinear-optimization. It's applicability to large-scale programming, however, is often limited due to the run-time complexity in computing the Newton step.

For a generic function \(\obj{f} : \mathbb{R}^N \rightarrow \mathbb{R}\), computing the Hessian requires atleast \(O(N^2)\) flops; further inverting the matrix requires \(O(N^{\gamma})\) flops (\(\gamma = 3\), in practice). This is computationally infeasible for many problems in practice. 

Often, however, one is also given access to the the \emph{computational structure} of the objective. The computer routine for calculating the objective \(\obj{f}(\cdot)\) can be represented as a Directed Acyclic Graph [DAG] mapping inputs to \(\obj{f}(\cdot)\) \emph{via} intermediary nodes. 

For instance, the objective function for the canonical optimal-control problem is given by,
\begin{equation}
\label{eqn:optcon}
\begin{aligned}
&\min_{u_0, u_1, \dots, u_{n - 1}} \left[\obj{J}(u_0, \dots, u_n) \triangleq \sum_{i = 0}^{n - 1} \obj{l}_i(x_i, u_i) + \obj{l}_n(x_n)\right],\\
&\quad \forall i, x_{i + 1} \leftarrow \mathtt{f}(x_i, u_i),
\end{aligned}
\end{equation}
where the dynamics and local-objectives of the system are given by \(\mathtt{f}(\cdot, \cdot)\), and \(\obj{l}_i(\cdot, \cdot)\) respectively. The infix operator '\(\leftarrow\)' indicates that the value appearing on the right-hand side, is given the placeholder symbol present to its left; we explicitly distinguish this from the '\(=\)' operator, which is taken to represent a constraint.

\begin{figure}[h!]
  \label{fig:optcon-directed}
  \centering
  \begin{tikzpicture}[shorten >=1pt,->]
    \tikzstyle{vertex}=[circle, minimum size=0pt,inner sep=0pt]

    \foreach \name/\x in {0/0, 1/1, 2/2, 3/3, n-1/6}
    \node[vertex] (U-\name) at (\x,0) {$u_{\name}$};   

    \foreach \name/\x in {1/1, 2/2, 3/3, n-1/6, n/7}
    \node[vertex] (X-\name) at (\x, 25pt) {$x_{\name}$};

    % \node[vertex] (V) at (8, 25pt) {$V_{n}$};
    
    \foreach \from/\to in {1/2, 2/3, n-1/n}
    { \draw (X-\from) edge (X-\to); }
    % \draw (X-n) edge (V);
    \draw[->, dotted] (X-3) -- (4, 25pt);
    \draw[->, dotted] (X-3) -- (4, 25pt);
    \draw[->, dotted] (5, 25pt) -- (X-n-1);

    \foreach \from/\to in {0/1, 1/2, 2/3, n-1/n}
    { \draw[->] (U-\from) -- (X-\to); }
    \draw[-, dotted] (U-3) -- (4, 25pt);
    \draw[-, dotted] (4, 0) -- (5, 0);
  \end{tikzpicture}
\caption{Optimal control problem: The dynamical system states are represented by \(\{x_i\}\), and the control by nodes \(\{u_i\}\).}
\end{figure}
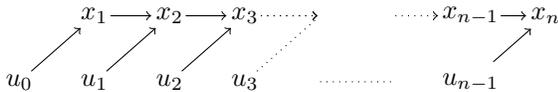
The order of computation for the objective \eqnref{eqn:optcon} can be represented by a linear-chain \figref{fig:optcon-directed}. Lacking constraints, the apparent sparsity in \eqnref{eqn:optcon}, is entirely destroyed once all the placeholders are substituted for,
\begin{flalign*}
\obj{J}(u_0, \dots, u_n) = &\obj{l}_0(x_0, u_0) +\\
&\obj{l}_1(\mathtt{f}(x_0, u_0), u_1) +\\
&\obj{l}_2(\mathtt{f}(\mathtt{f}(x_0, u_0), u_1), u_2) + \dots.
\end{flalign*}
The Hessian of \(\obj{J}(\cdot)\) thus being dense, implies a run-time that is cubic in the input dimensions for the Newton step computation; computing the Hessian itself is quadratic.

By contrast, once the problem \eqnref{eqn:optcon} is written in its constrained form (by replacing '\(\leftarrow\)' with '\(=\)'), the sparsity of the resulting Karush-Kuhn-Tucker [KKT] system, readily allows for computing the SQP/Lagrange-Newton step in linear time \cite{wright1990solution}. Such a transformation, however, comes at the cost of increasing the size of the optimization problem, abandoning state feasibility, and increased implementation complexity.

The question which this paper answers, is whether there exist \emph{general} techniques, which allow exploiting the sparsity of the problem, while working solely with the input variables. Note that these are not questions merely about elimination orders, but are also verily algebraic in nature.

\textsc{Automatic Differentiation}: Research on Automatic Differentiation [AD] has produced many techniques for exploiting the \emph{computational structure} of generic functions. They are routinely employed for efficient calculation of gradients and Hessian vector products \cite{griewank2008evaluating}. The applicability of AD to second-order optimization is, however, quite limited. 

AD is typically used either for computing the entire Hessian matrix, or for calculating Hessian vector products for use in Nonlinear Conjugate Gradient Descent [CG]. Hessians are computed, by accumulating one column at a time \emph{via} calls to the Hessian vector product routine \cite{griewank2008evaluating}. The sparsity of the Hessian can be exploited in reducing the number of such calls \cite{coleman1984estimation}, but structured problems such as \eqnref{eqn:optcon} will not allow for any such economy. Compositional chains of functions, such as those in optimal control \eqnref{eqn:optcon} \figref{fig:optcon-directed}, not only serve to make Hessians dense but can also lead to condition numbers, exponential in their diameters. Large condition numbers are likely to negate any computational advantages offered by methods like CG. 

The above techniques form the Hessian matrix, directly or indirectly, before computing the Newton step. This stands in contrast with the root-finding problem, for which there do exist methods for directly computing the Newton-Raphson step \cite{utke1996efficient} \cite{griewank2008evaluating} \cite{dixon2009automatic}. The root finding problem involves the inversion of the Jacobian of a function - rather than the Hessian - and these methods reduce this computation to that of inverting a sparse matrix \cite{utke1996efficient}. The Newton step (for optimization) can also be computed using this method by formulating it as a root-finding problem on the gradient. This, however, results in non-symmetric matrices depending on the computational graph of the gradient, as opposed to the function itself. The latter graph is transitively closed, and hence, analyses for the above Newton-Raphson AD algorithm apply to the gradient, and are difficult to extend to the underlying objective \cite{dixon2009automatic}.

\textsc{Dynamic Programming}: The question posed earlier, has already been answered in the affirmitive, for the optimal control problem. There exists an algorithm for optimal control, based on Dynamic Programming, that exploits algebraic dependencies in \eqnref{eqn:optcon}, in order to compute the Newton-step in only linear time \cite{de1988differential} \cite{jacobson1970differential} \cite{todorov2005generalized}. 

The run-time of this algorithm is the direct result of the sparsity of the corresponding constrained problem  \cite{de1988differential} \cite{wright1990solution}\cite{ralph1996parallel}. The band-structure of the relevant KKT system allows for solving the system in linear time \cite{wright1990solution}. The relationship between computing the Newton-step (Hessian of the objective), and computing the Lagrange-Newton step (Hessian of the Lagrangian), is established by noting that there exist multiplier values such that both compute the \emph{same} result \cite{de1988differential}.

Such algorithms are routinely employed by practitioners for updating \emph{control policies} in real-time, while maintaining a feasible trajectory. These algorithms have been extended to Extended Kalman Filtering (EKF), as well as various other formulations of the control problem \cite{toussaint2010bayesian} \cite{tassa2011optimal}  \cite{wright1993interior}.

%\cite{bell2009inequality}
% Computationally, each call to the Hessian vector product AD routine, is nearly as expensive as computing the Newton step through such algorithms \cite{griewank2008evaluating} \cite{christianson1999cheap}. Hence, these are paradigmatically different from current AD approaches.

\textsc{Overview}: We generalize such algorithms, by using Hessian vector product equations from AD, to relate the computation of Newton step and Lagrange-Newton step, for arbitrary structured objectives. 

We then extend this framework to structured optimization problems with equality constraints.

Further, we show that solving the resultant KKT systems can be accomplished in time \(\tilde{O}(\tw^3)\), where '\(\tw\)' is the tree-width of the canonical computational graph.

Finally, we show results from numerical experiments.

\section{Notation}
Let \(\mathcal{G}\) be a Directed Acyclic Graph [DAG], and let each vertex \(v \in V[\mathcal{G}]\), be associated with \emph{state} \(\pr{v} \in U_v \subset \mathbb{R}^{n_v}\), taking values in an open set. Denote by \(\pa(v)\), the parents of \(v \in V[\mathcal{G}]\), and by \(\ch(v)\) its children; let \(\pr{A}\) be the (labelled) concatenation of \emph{state}s, associated with vertices in set \(A \subset V[\mathcal{G}]\). Define the set of \emph{input} nodes \(X = \{x_1, x_2,\dots, x_n\} \triangleq \{v \mid \pa(v) = \emptyset, v \in V[\mathcal{G}]\}\), to be the parentless vertices of \(\mathcal{G}\).

An objective function \(\obj{f} : U_{x_1} \times \dots U_{x_n} \rightarrow \mathbb{R}\), has the \emph{computational structure} given by the tuple \((\mathcal{G}, \{\varphi_v\}, \{\obj{l}_v\})\), if it can be written as the sum of local objectives \(\obj{l}_v : \prod_{z \in \{v\} \cup \pa(v)} U_z \rightarrow \mathbb{R}\), on the graph \(\mathcal{G}\),
\begin{equation}
\label{eqn:struct}
\begin{aligned}
\obj{f}: (\pr{x_1}, \dots, \pr{x_n})& \mapsto \sum_{v \in V[\mathcal{G}]} \obj{l}_v(\pr{v \cup \pa(v)}),\\
\pr{v} \leftarrow \varphi_v(\pr{\pa(v)})&, \quad \forall v \in V[\mathcal{G}], \pa(v) \neq \emptyset.
\end{aligned}
\end{equation}
The \emph{state} of a non-input node \(v \in V[\mathcal{G}]\) in \eqnref{eqn:struct}, is defined recursively as \(\pr{v} \leftarrow \varphi_v(\pr{\pa(v)})\), for some given function \(\varphi_v : \prod_{z \in \pa(v)} U_z \rightarrow U_v\). It follows since \(\mathcal{G}\) is a DAG, that \(\pr{V[\mathcal{G}]}\) and hence \(\obj{f}(\cdot)\), is uniquely determined from the input \(\pr{X}\), and functions \(\{\varphi_v\}\). The order of computation for the objective is given by the topological ordering of \(\mathcal{G}\), and the DAG \(\mathcal{G}\) is called the \emph{computational graph} of \(\obj{f}(\cdot)\). The computer routine for calculating any objective function, can be represented by such a structure \cite{griewank2008evaluating}.

In the following sections, the symbolism \(\p_u v\) is used as a shorthand for \({\p \pr{v} \over \p \pr{u}} \big|_{\pr{X}}\). The derivatives of functions with respect to \(\pr{u}\) are similarly denoted by the operator \(\p_{u}\); that with respect to a (labelled) set \(A = \{v_1, v_2, \dots\} \subset V[\mathcal{G}]\) by \(\p_{A} \triangleq [\p_{a_1}, \p_{a_2}, \dots]\).

\section{Newton step}
\label{sec:newt}
Consider the objective function in \eqnref{eqn:struct}, defined by the tuple \((\mathcal{G}, \{\varphi_v\}, \{\obj{l}_v\})\). The optimization problem of interest is the following,
\begin{equation}
  \label{eqn:obj}
\begin{aligned}
\min_{\pr{x_1}, \dots, \pr{x_n}} &\left( \obj{f} \triangleq \sum_{v \in V[\mathcal{G}]} \obj{l}_v(\pr{v \cup \pa(v)}) \right),\\
\pr{v} \leftarrow \varphi_v(\pr{\pa(v)}),& \quad \forall v \in V[\mathcal{G}], \pa(v) \neq \emptyset,
\end{aligned}
\end{equation}
and the corresponding constrained problem is obtained by replacing the operator '\(\leftarrow\)' by '\(=\)' in \eqnref{eqn:obj}.

In the following, we consider first the constrained formulation of \eqnref{eqn:obj}, and define the KKT system involved in computing the Lagrange-Newton step; we then relate these to computing the Newton step.

% \begin{equation}
%   \label{eqn:cobj}
% \begin{aligned}
% \min_{\pr{V}} &\left(\sum_{v \in V[\mathcal{G}]} \obj{l}_v(\pr{v \cup \pa(v)}) \right),\\
% \varphi_v(\pr{\pa(v)}) - &\pr{v} = 0, \quad \forall v \in V[\mathcal{G}], \pa(v) \neq \emptyset.
% \end{aligned}
% \end{equation}
%We assume that the Hessian of \(\obj{f}\) is positive-definite [p.d]. 

%\footnote{We note that strictly speaking, \(\p_u f\) does not make sense for any \(u\) that is not an input to the graph \(\mathcal{G}\). In such a case it should be interpreted as cutting off vertex \(u\) from it ancestors, and elevating it to the status of an independent variable, for the modified function.}

\subsection{Lagrange-Newton}
The Lagrangian for the constrained form of (\ref{eqn:obj}), is given by,
\begin{equation}
  \label{eqn:cdef}
  \begin{aligned}
    &\mathcal{L}(\pr{V[\mathcal{G}]}, \lambda) \triangleq\\ 
    &\quad\quad\sum_{v \in V[\mathcal{G}]} \obj{l}_v(\pr{v \cup \pa(v)}) + \sum_{\substack{v \in V[\mathcal{G}],\\ \pa(v) \neq \emptyset}} \lambda_v^{\mathtt{T}} h_v(\pr{v \cup \pa(v)}),\\
    &\mbox{where,} \\
    &\forall v \in V[\mathcal{G}], \pa(v) \neq \emptyset, \quad h_v(\pr{v \cup \pa(v)}) \triangleq \varphi_v(\pr{\pa(v)}) - \pr{v},
  \end{aligned}
\end{equation}
and the vector \(\lambda\) is the labelled concatenation of all \(\lambda_v\)'s.
% \begin{equation}
%   \min_{\pr(V(\mathcal{G}))} \mathcal{L}
%   \p_{V} \mathcal{L}(\pr{V}^*, \lambda^*) = 0,\quad g(\pr{V}^*) = 0.
% \end{equation}

The necessary first order conditions for optimality of this problem are given by \cite{nocedal2006numerical},
\begin{equation}
  \label{eqn:lagrangefirst}
  \p_{V} \mathcal{L}(\pr{V}^*, \lambda^*) = 0,\quad h(\pr{V}^*) = 0.
\end{equation}
The Lagrange-Newton step for solving this system of equations, around a nominal \((\tilde{\pr{V}}, \lambda)\), entails solving the following KKT system \cite{nocedal2006numerical},
% Taking a first order variations of (\ref{eqn:lagrangefirst}) w.r.t \((\delta \pr{V}, \delta \lambda)\), around a nominal \(\tilde{\pr{V}}\), results in the following KKT system,
\begin{equation}
  \label{eqn:clnstep}
  \left[\begin{array}{c c}
          \p^2_{V} \mathcal{L} & \p_{V} h^{\mathtt{T}} \\
          \p_{V} h & 0
    \end{array}\right] \left[\begin{array}{c} \delta \pr{V} \\ \delta \lambda
    \end{array}\right] = \left[\begin{array}{c} - \p_{V} \mathcal{L} \\ -h
    \end{array} \right].
\end{equation}
% \begin{equation}
%   \begin{aligned}
%     & \p^2_{V} \mathcal{L} \cdot \delta \pr{V} + \p_{V} g \cdot \delta \lambda = - \p_{V} \mathcal{L},\\
%     &\p_{V} g \cdot \delta \pr{V} = - g .
%   \end{aligned}
% \end{equation}
% These equations can be solved by introducing another Lagrange multiplier \(\delta \lambda\),

Sequential Quadratic Programming [SQP], involves taking a step along \((\delta \pr{V}, \delta \lambda)\) and iteratively solving for the first order conditions. In the following section, it will be shown that there exist values for Lagrange multipliers, depending only on the inputs, such that the solution to \eqnref{eqn:clnstep}, yields the Newton step for the unconstrained objective. 

%This obviates the need for storing \& updating multipliers, and provides the algebraic.

% It will be seen that if the Lagrange multipliers are deterministically set to \(\lambda_v = \p_v f, \forall v \in V[\mathcal{G}]\), the obtained Lagrange-Newton step from (\ref{eqn:clnstep}) yields the Newton step for the original objective in (\ref{eqn:obj}). This relationship between SQP multipliers and the descent direction was first noticed for the trajectory optimization problem in Ralph's seminal paper \cite{ralph1996parallel}.

% More generally, the Lagrange-Newton step for every constrained optimization problem which factorizes according to \(\mathcal{H}\),
% \begin{equation}
%   \label{eqn:cobj}
%   \min_{\pr{X}} \sum_{c \in E[\mathcal{H}]} \obj{l}_c(\pr{c}), \quad g_{c}(\pr{c}) = 0, \forall c \in E[\mathcal{H}],
% \end{equation}
% can also be computed in \(\tilde{O}(\tw(\mathcal{H})^3)\) time. This is interesting in that it reflects here the dichotomy that exists between Bayesian \& Markov networks in graphical models.
\subsection{Unconstrained Newton}
We recollect certain defintions from AD, and then continue to present one of the central results of the paper.

\textsc{Reverse AD}: The first derivatives of the objective \(\obj{f}(\cdot)\) can be calculated by applying the chain rule over \(\mathcal{G}\),
\begin{equation}
  \label{eqn:rmad}
  \begin{aligned}
    %& \p_v f = \p_v \obj{l}_v + \sum_{k: \arcp{v}{k}} (\p_k f \cdot \p_vk + \p_v \obj{l}_k); \quad \arc{v}{k} \in E[\mathcal{G}] \Rightarrow \p_v k \triangleq {\p \varphi_k(\pr{\pa(k)}) \over \p \pr{v}}.\\
    &\forall v,\quad  \p_v \obj{f} = \sum_{s \in v \cup \ch(v)} \p_v \obj{l}_s + \sum_{d \in \ch(v)} \p_d \obj{f}^{\mathtt{T}} \; \p_v d;\\
    &\quad\quad\quad v \in \pa(d) \Rightarrow \p_v d \triangleq {\p \varphi_d(\pr{\pa(d)}) \over \p \pr{v}}.
  \end{aligned}
\end{equation}
Since \(\mathcal{G}\) is a DAG, there exist child-less nodes (\emph{i.e} \(\ch(v) = \emptyset\)), from which the above recursion can be initialized. The recursion then proceeds backward in the depth first search order on \(\mathcal{G}\). This algorithm is known as reverse-mode AD \cite{griewank2008evaluating}.

% \begin{equation}
%   \label{eqn:rmad}
%   \begin{aligned}
%     &\mbox{Backward recursion}:\\
%     %& \p_v f = \p_v \obj{l}_v + \sum_{k: \arcp{v}{k}} (\p_k f \cdot \p_vk + \p_v \obj{l}_k); \quad \arc{v}{k} \in E[\mathcal{G}] \Rightarrow \p_v k \triangleq {\p \varphi_k(\pr{\pa(k)}) \over \p \pr{v}}.\\
%     & \p_v \obj{f} = \sum_{s \in v \cup \ch(v)} \p_v \obj{l}_s + \sum_{d \in \ch(v)} \p_d \obj{f} \cdot \p_v d; \quad \arc{v}{w} \in E[\mathcal{G}] \Rightarrow \p_v w \triangleq {\p \varphi_w(\pr{\pa(w)}) \over \p \pr{v}}.
%   \end{aligned}
% \end{equation}
\textsc{Hessian vector AD}: A change in the inputs \(\delta {\pr{X}}\), results in the first-order change in the derivative, \(\delta [\p_v \obj{f}] \triangleq \p_{X v}^2 \obj{f} \cdot \delta \pr{X}\), which is given by the Hessian vector product. Computing the Newton step is thus, equivalent to finding a \(\delta {\pr{X}}\) such that, \(\delta [\p_X \obj{f}] = - \p_X \obj{f}\). 

Applying chain-rule over the DAG \(\mathcal{G}\), for all terms in \eqnref{eqn:rmad}, we obtain,
\begin{equation}
  \label{eqn:armad}
  \begin{aligned}
    % &\mbox{Backward recursion}:\\
    &\forall v,\quad  \delta [\p_v \obj{f}] =\\&\quad\quad \sum_{s \in v \cup \ch(v)} \left(\sum_{a \in v \cup \pa(s)} \p_{av}^2 \obj{l}_s \cdot \delta \pr{a} \right) +\\
    &\quad\quad \sum_{d \in \ch(v)} \left(\delta [\p_d \obj{f}]^{\mathtt{T}} \p_v d + \sum_{a \in \pa(d)} (\p_d \obj{f}^{\mathtt{T}} \; \p_{av}^2 d) \cdot \delta \pr{a})\right);\\
    %+ \sum_{k: \arcp{v}{k}} \delta [\p_k f] \cdot \p_v k\; + \\
    %\sum_{k: \arcp{v}{k}} &\left( \sum_{w \in k \cup \pa(k)} \p_{vw} \obj{l}_k \cdot \delta \pr{w} + \sum_{w \in \pa(k)} (\p_k f \cdot \p_{v w} k) \cdot \delta \pr{w} \right),\\
    %&\mbox{Forward recursion}:\\    
    &\forall a,\quad \delta \pr{a} = \sum_{d \in \pa(a)} \p_{d} a \cdot \delta \pr{{d}}.
  \end{aligned}
\end{equation}
These equations can be solved, for a given \(\delta \pr{X}\), by a forward-backward recursion similar to the one used for solving \eqnref{eqn:rmad} \cite{griewank2008evaluating}. Computing the Hessian-vector product in this manner takes time \(\tilde{O}(\omega(\mathcal{\hat{G}})^2)\) \footnote{We use \(\tilde{O}(.)\) to hide factors linear in \(|E| + |V|\).} \cite{griewank2008evaluating}, where \(\omega(\mathcal{\hat{G}})\) is the clique number of the moralization of \(\mathcal{G}\).
% \footnote{The clique number \(\omega(\mathcal{G})\) of a direc equals that of its primal graph \(\mathcal{G}_p \triangleq (V[\mathcal{H}], \{uv | \exists e \in E[\mathcal{H}], u, v \in e\}).\)}

\textsc{Newton step}: The problem of interest is, however, the exact inverse: find a \(\delta \pr{X}\), such that \(\delta[\p_X \obj{f}] = -\p_X \obj{f}\). This question is answered by the following theorem.
\begin{theo}
\label{theo:gnewt}(Newton step) The Newton step for the objective \eqnref{eqn:struct} is given by the Lagrange-Newton step \eqnref{eqn:clnstep}, when \(\pr{V}\) is feasible and when \(\forall v, \lambda_v = \p_v \obj{f}\) as defined in \eqnref{eqn:rmad}.
\end{theo}
\emph{Proof.}
The second equation in \eqnref{eqn:armad} is equivalent to \(\p_V h \cdot \delta S_V = -h\), in \eqnref{eqn:clnstep}. Rearranging the first equation from \eqnref{eqn:armad}, and setting \(\delta[\p_v \obj{f}] = -\p_v \obj{f}\) for all inputs, we obtain \(\forall v\),
\begin{equation}
\begin{aligned}
\label{eqn:ppobj}
0 =&\\
&\sum_{\substack{s \in v \cup \ch(v),\\ a \in v \cup \pa(s)}} \p_{va}^2 \obj{l}_s \;\delta \pr{a} + \sum_{\substack{d \in \ch(v),\\ a \in \pa(d)}} (\p_d \obj{f}^{\mathtt{T}}\; \p_{av}^2 d) \; \delta \pr{a} +\\
& - \left(\begin{cases} \delta [\p_v \obj{f}] & \pa(v) \neq \emptyset\\ - \p_v \obj{f} & \mbox{otherwise}\end{cases}\right) + \sum_{d \in \ch(v)}  (\p_v d)^{\mathtt{T}} \delta [\p_d \obj{f}].
\end{aligned}
\end{equation}
Similarly, expanding the top block in \eqnref{eqn:clnstep} using the definitions in \eqnref{eqn:obj} \& \eqnref{eqn:clnstep}, we obtain \(\forall v\), 
\begin{equation}
\begin{aligned}
\label{eqn:ppLexp}
&-\p_{v} \mathcal{L} =\\
&\sum_{\substack{s \in v \cup \ch(v),\\a \in v \cup \pa(s)}} \p_{va}^2 \obj{l}_s \; \delta \pr{a} + \sum_{\substack{d \in \ch(v),\\a \in \pa(d)}} (\lambda_d^{\mathtt{T}} \; \p_{av}^2 d) \; \delta \pr{a} +\\
& - \left(\begin{cases} \delta \lambda_v & \pa(v) \neq \emptyset \\ 0 & \mbox{otherwise}\end{cases} \right) + \sum_{d \in \ch(v)}  (\p_v d)^{\mathtt{T}} \; \delta \lambda_d ,
%\sum_{s \in v \cup \ch(v)}  \p_v \obj{l}_s + \sum_{d \in \ch(v)} \lambda_d \cdot \p_{v} d = \p_{v} \obj{f}
% + \sum_{a \in \pa(w)} (\p_d \obj{f} \cdot \p_{av}^2 d) \cdot \delta \pr{a})\right)\\
\end{aligned}
\end{equation}
where,
\begin{equation}
\begin{aligned}
\label{eqn:pLexp}
&\p_{v} \mathcal{L} = \\
&\begin{cases} \sum_{s \in v \cup \ch(v)}  \p_v \obj{l}_s + \sum_{d \in \ch(v)} \lambda_d \cdot \p_{v} d, & \pa(v) = \emptyset\\ \begin{array}{l} \sum_{s \in v \cup \ch(v)}  \p_v \obj{l}_s + \sum_{d \in \ch(v)} \lambda_d \cdot \p_{v} d \\- \lambda_v \end{array}, & \mbox{otherwise} \\  \end{cases}
\end{aligned}
\end{equation}
% \textcolor{red}{\textbf{---------------------------------------------------------}}
The result follows from equations \eqnref{eqn:rmad}, \eqnref{eqn:ppobj}, \eqnref{eqn:ppLexp} \& (\ref{eqn:pLexp}).

\qed

% We also see that the terms \(\delta[\p_v \obj{f}]\) in (\ref{eqn:ppobj}) are equivalent to \(\delta \lambda_v\) in (\ref{eqn:clnstep}), when the assumptions of Theorem~\ref{theo:gnewt} hold.

\textsc{Graphical Newton}: The above theorem immediately yields the following optimization algorithm,
\begin{algorithm}[h]
  \caption{Graphical Newton}
  \label{alg:gnewt}
\begin{algorithmic}[1]
  \STATE {\bfseries Input:} initial $\pr{X}$, tuple $(\mathcal{G}, \{\varphi_v\}, \{\obj{l}_v\})$
  \REPEAT 
  \STATE Compute \(\obj{f}, \{\p_v \obj{f}\}, \{\p^2 \varphi_v\}\) from (\ref{eqn:struct}), (\ref{eqn:rmad}).
  \STATE Compute the SQP step from (\ref{eqn:clnstep}), with \(\lambda_v = \p_v \obj{f}, \forall v\).
  \STATE Compute step-length \(\eta\) \emph{via} linesearch on inputs \(\pr{X}\).
  \STATE Update inputs: $\pr{X} \leftarrow \pr{X} + \eta \delta \pr{X}$.
  \UNTIL{$\norm{\p_X \obj{f}} \leq \epsilon$ }
\end{algorithmic}
\end{algorithm}

%%FIX FIX FIX
The run-time of every iteration in Algorithm~\ref{alg:gnewt} depends crucially upon the time required to solve (\ref{eqn:clnstep}). The run-time bounds for solving such KKT systems is taken up later in the paper.

% Complexities of solving sparse linear systems are often given in terms of the \emph{tree-width} of the support of the matrix.

 % such theorems are not applicable when the diagonal is zero \cite{bridson2007ordering} \cite{davis2006direct}. In the following section, we show that (\ref{eqn:clnstep}) can be solved in \(\tilde{O}(\tw(\mathcal{G})^3)\).

%The Newton search direction is given by the input \(\delta \pr{X}^{\dagger}\) that satisfies \(\delta [\p_X \obj{f}] \triangleq [\p_{XX}^2 \obj{f} ] \cdot \delta \pr{X}^{\dagger} = -\p_{V} f\). If the objective function \(\obj{f}(\cdot)\) were quadratic, then the iterate \(\pr{X}^+ \triangleq \pr{X} + \delta \pr{X}^{\dagger}\) would be its global minimizer. Otherwise, this iteration scheme converges to a local minimum, quadratically \cite{nocedal2006numerical}.

\subsection{Extension to equality constraints}
Consider optimization problems, which have equality constraints in addition to the structured objective from before,
\begin{equation}
  \label{eqn:cobj}
\begin{aligned}
\min_{\pr{x_1}, \dots, \pr{x_n}} &\left( \obj{f} \triangleq \sum_{v \in V[\mathcal{G}]} \obj{l}_v(\pr{v \cup \pa(v)}) \right),\\
\pr{v} \leftarrow \varphi_v(\pr{\pa(v)}),& \quad \forall v \in V[\mathcal{G}], \pa(v) \neq \emptyset,\\
& c(\pr{C}) = 0,
\end{aligned}
\end{equation}
where \(c(\cdot) = 0\) is an additional equality constraint, which depends on the variables \(C \subset V[\mathcal{G}]\). The Lagrangian for this problem is given by,
\begin{equation}
  \label{eqn:cdef}
  \hat{\mathcal{L}}(\pr{V[\mathcal{G}]}, \lambda) = \mathcal{L}(\pr{V[\mathcal{G}]}, \lambda_{V \backslash X}) + \lambda_h^{\mathtt{T}} c(\pr{C}),
\end{equation}
where \(\mathcal{L}\) is as defined in \eqnref{eqn:cdef}, and \(\lambda_{V \backslash X}\) is the corresponding set of multipliers; the variable \(\lambda\), being the concatenation of \(\lambda_c\) and all multipliers, \(\lambda_{V\backslash X}\), appearing in \eqnref{eqn:cdef}.

Theorem~\ref{theo:gnewt} can be applied to this problem by treating \(\lambda_c^T c(\pr{C})\) as another cost function in the objective, while also including the constraint in the KKT system \eqnref{eqn:clnstep}. The iteration can then proceed by solving the KKT system with \(\lambda_v = \p_v (\obj{f} + \lambda_c^T c), \forall v\), and using a merit function for the linesearch procedure; the variables \((\pr{X}, \lambda_c)\) are updated accordingly. We omit the proof for the validity of this method.

\section{Message Passing}
The classical run-time bound for Cholesky factorization (\emph{i.e} Gaussian Belief Propagation \footnote{Gaussian-BP, computes the LU decomposition of a matrix}) \cite{davis2006direct} \cite{wainwright2008graphical}, cannot be extended to problems such as \eqnref{eqn:clnstep}, because of the appearance of linear constraints. Such bounds for structured KKT systems, do not appear to be known within the sparse linear algebra community \cite{bridson2007ordering}.

In this section, we provide a Message Passing algorithm for solving such KKT systems, and show that it has a run-time bound of \(\tilde{O}(\tw^3)\) \footnote{The tilde hides factors linear in \(|V[\mathcal{H}]|, |E[\mathcal{H}]|\).}, given the tree-decomposition.

% We note that since LDL decomposition, unlike Cholesky, is not numerically stable, Message Passing schemes such as those employed in \cite{toussaint2010bayesian} may be ill-conditioned.

%Sophisticated packages such as MA57 \cite{duff2004ma57} incorporate a variety of techniques in order to ensure numerical stability. While theoretically, LDL 

\subsection{Hypergraph structured QPs}
For a hypergraph \(\mathcal{H}\), denote the adjacency and incidence matrices by \(\mathcal{A}[\mathcal{H}]\) \& \(\mathcal{B}[\mathcal{H}]\) respectively,
\begin{equation}
\begin{aligned}
  &\mathcal{A}[\mathcal{H}] \in \mathbb{R}^{|V[\mathcal{H}]|\times|V[\mathcal{H}]|}, \quad\quad \mathcal{B}[\mathcal{H}] \in \mathbb{R}^{|E[\mathcal{H}]|\times|V[\mathcal{H}]|},\\
  &\mathcal{A}[\mathcal{H}]_{uv} = \begin{cases} 1 & \exists e \in E[\mathcal{H}],\; u, v \in e\\ 0 & \mbox{otherwise}\end{cases}\\  
  &\mathcal{B}[\mathcal{H}]_{eu} = \begin{cases} 1 & u \in e\\ 0 & \mbox{otherwise}\end{cases}
\end{aligned}
\end{equation}
Given such a hypergraph \(\mathcal{H}\), the family of QPs we're interested in solving is the following,
\begin{equation}
  \label{eqn:spqpfamily}
  \begin{aligned}
    &\min_{x} \sum_{e \in E[\mathcal{H}]}{1 \over 2} \pr{e}^{\mathtt{T}} Q_e \pr{e} - b_e^{\mathtt{T}} \pr{e}, \\ 
    &\forall e \in E[\mathcal{H}], \quad G_{e} \pr{e} = h_e.
  \end{aligned}
\end{equation}
Assuming that the QP has a bounded solution and that the constraints are full rank, the minimizer to (\ref{eqn:spqpfamily}) is given by the solution to the following KKT system,
\begin{equation}
  \label{eqn:spfamily}
  \begin{aligned}
    \left[\begin{array}{c c} Q & G^{\mathtt{T}} \\ G & 0\end{array}\right] \left[\begin{array}{c} x \\ \lambda\end{array}\right] &= 
    \left[\begin{array}{c} b \\ h \end{array}\right],\\  
    x, b \in \mathbb{R}^{|V|}, \quad& \lambda, h \in \mathbb{R}^{M},\\
\end{aligned}
\end{equation}
where \(Q, G, \lambda, x, b\) are concatenation of terms defined in \eqnref{eqn:spqpfamily} respectively. The sparsity/support of (\ref{eqn:spfamily}) is closely related to \(\mathcal{H}\), since,
\begin{flalign*}
&\supp(Q) \subseteq \supp(\mathcal{A}[\mathcal{H}]),\\
&\forall i, \exists e, \supp(G_{i, :}) \subseteq \supp(\mathcal{B}[\mathcal{H}]_{e, :}).
\end{flalign*}
Every row of the constraint, \(G_{i, :}\), has the same sparsity as some edge \(e \in E[\mathcal{H}]\). %We note that it may be easier to read this section with the notion of undirected graphs, and their corresponding set of maximal cliques.

\textsc{Tree decomposition}:
%\subsubsection{Tree decomposition}
Extending the notion of Dynamic Programming to non-trees (including Hypergraphs) requires a partitioning of the graph so as to satisfy a \emph{lifted} notion of being a tree \cite{kleinberg2006algorithm}. Tree decomposition captures the essence of such graph partitions,

\begin{defn}{(Tree decomposition)}
A tree-decomposition of a hypergraph \(\mathcal{H}\) consists of a tree \(\mathcal{T}\) and a map \(\chi : V[\mathcal{T}] \rightarrow 2^{V[\mathcal{H}]}\), such that,
\begin{enumerate}[i]
\item (\emph{Vertex cover}) \(\cup_{i \in V[\mathcal{T}]} \chi(i) = V[\mathcal{H}].\)
\item (\emph{Edge cover}) \(\forall e \in E[\mathcal{H}],\; \exists i \in V[\mathcal{T}], e \subset \chi(i).\)
\item (\emph{Induced sub-tree}) \(\forall u \in V[\mathcal{H}],\; \mathcal{T}_u \triangleq \mathcal{T}[\{i \in V[\mathcal{T}]| u \in \chi(i)\}]\; \mbox{is a non-empty subtree}\)
\end{enumerate}
The tree-width of a tree-decomposition \(\mathcal{T}\) is defined to be \(\tw(\mathcal{T}) = \max_{v \in V[\mathcal{T}]} |\chi(v)| - 1\). The tree-width of a graph \(\mathcal{H}\) is defined to be the minimal tree-width attained by any tree-decomposition of \(\mathcal{H}\).
\end{defn}

We define the vertex-induced subgraph in what follows to be \(\mathcal{H}[S] \triangleq (V[\mathcal{H}], \{e \cap S, e\in E[\mathcal{H}]\})\). The following lemma ensures that such a decomposition ensures \emph{local dependence} \cite{kleinberg2006algorithm}.
\begin{lem}
\label{lem:esep}{(Edge separation)} Deleting the edge \(xy \in E[\mathcal{T}]\), renders \(\mathcal{H}[V \backslash (\chi(x) \cap \chi(y))]\) disconnected.
\end{lem}
% It follows from the edge-separation theorem (\ref{theo:esep}), that \emph{Dynamic programming} algorithms on trees, can be lifted to tree-decompositions of general graphs.

\textsc{Hypertree structured QP}: The tree-decomposition itself can be considered a Hypergraph, \((V[\mathcal{H}], \{\chi(u), \forall u \in V[\mathcal{T}]\})\). Such a \emph{Hypertree}\footnote{There are multiple definitions of a \emph{Hypertree}; we use the term to mean a maximal Hypergraph, whose tree-decomposition can be expressed in terms of its edges.} can also be thought of as a Chordal graph \cite{wainwright2008graphical}. We assume henceforth that the given graph \(\mathcal{H}\) is a hypertree, and that \(\mathcal{T}\) is its tree-decomposition. 

\begin{algorithm}[h]
  \caption{Graphical QP}
  \label{alg:gqp}
\begin{algorithmic}[1]
  \STATE {\bfseries Given}: $\mathcal{T}, \mathcal{H}, \{Q_e\}, \{b_e\}, \{G_e\}, \{h_e\}.$
  \STATE
  \STATE {\bfseries function} GatherMessage($l, p, \mathcal{T}$)
  \STATE $(\tilde{Q}_l, \tilde{b}_l, \tilde{G}_l, \tilde{h}_l) \leftarrow (Q_l, b_l, G_l, h_l)$
  \FOR{$c \in \nh_{\mathcal{T}}(l)\backslash p$}
  \STATE $(Q_{c\rightarrow l}, G_{c\rightarrow l}, b_{c\rightarrow l}, h_{c\rightarrow l}) \leftarrow \mbox{GatherMessage}(c, p, \mathcal{T})$
  \STATE $(\tilde{Q_l}, \tilde{b_l}) \leftarrow (\tilde{Q_l}, \tilde{b_l}) + (Q_{c\rightarrow l}, b_{c\rightarrow l})$
  \STATE $\tilde{G}_l \leftarrow [\tilde{G}_l; G_{c\rightarrow l}], \tilde{h}_l \leftarrow [\tilde{h}_l; h_{c\rightarrow l}]$
  \ENDFOR
  \STATE {\bfseries return} Factorize($\chi(l), \chi(p), \tilde{Q}_l, \tilde{b}_l, \tilde{G}_l, \tilde{h}_l$)
  \STATE
  \STATE {\bfseries function} Factorize($\chi(l), \chi(p), \tilde{Q}, \tilde{b}, \tilde{G}, \tilde{h}$)
  \STATE $(\xi, \iota) \leftarrow (\chi(l) \backslash \chi(p), \chi(l) \cap \chi(p))$
  \STATE $r \leftarrow \rank(\tilde{Q}_{\iota, \iota})$
  \STATE {\bfseries return} Gaussian-BP messages from \eqnref{eqn:mpcore}.
% {\bfseries return}
  % \ELSE
  % \STATE {\bfseries return}
  % \ENDIF
  \STATE {\bfseries return}
\end{algorithmic}
\end{algorithm}

% Add proof.
% \begin{lem}
% \label{lem:twlemqp}
% The QP (\ref{eqn:spqpfamily}) can be solved in time \(\tilde{O}(\tw(\mathcal{H})^3)\), given the minimal tree-decomposition.
% %\proof See the supplement.
% \end{lem}
% \proof Let \(\mathcal{H}\) be a Hypergraph as in (\ref{eqn:spqpfamily}), and let \((\mathcal{T}, \chi)\) be its tree-decomposition. We note that for every \(h \in E[\mathcal{H}]\), the number of constraints cannot exceed the tree-width.

%  Consider the leaf node \(l \in V[\mathcal{T}]\), and let \(p\) be its parent node in the tree. Consider the principal minor corresponding to \(p \& l\).

The gather stage of the Message Passing algorithm, is illustrated in Algorithm~\ref{alg:gqp}. \footnote{Note that the addition is performed vertex label-wise in Line~6 of Algorithm~\ref{alg:gqp}.} 

The function, Factorize, computes the partial LU decomposition of its arguments; we describe below, its operation. Denote the vertices that are interior to \(l\) by \(\iota = \chi(l) \cap \chi(p)\), and those on the boundary (\emph{i.e} common to \(p, l\)) by \(\xi = \chi(l) \backslash \chi(p)\), and let \(r = \rank(\tilde{Q}_{\iota, \iota})\).  The function computes Gaussian-BP messages from block pivots \(\{2, 3\}\) to \(\{1, 4\}\) in \eqnref{eqn:mpcore}. Note that, unlike Gaussian-BP, the matrices in \eqnref{eqn:mpcore} are not necessarily positive definite, but are however invertible.
% \begin{equation}
% \left[\begin{array}{c c c}
%         \tilde{Q}_{\xi \xi} & \tilde{Q}_{\iota \xi}^{\mathtt{T}} & \tilde{G}_{\xi}^{\mathtt{T}}\\
%         \tilde{Q}_{\iota \xi} & \tilde{Q}_{\iota \iota} & \tilde{G}_{\iota}^{\mathtt{T}}\\
%         \tilde{G}_{\xi} & \tilde{G}_{\iota} & 0
%       \end{array}\right] \left[\begin{array}{c}
%                                  x_l\\
%                                  x_p\\
%                                  \lambda
%                                \end{array}\right]  = \left[\begin{array}{c}
%                                                              \tilde{b}_{\xi}\\
%                                                              \tilde{b}_{\iota}\\
%                                                              \tilde{h}
%                                                            \end{array}\right]
% \end{equation}
\begin{figure*}[h!]
\label{fig:sd}
\centering
\includegraphics[height=5cm]{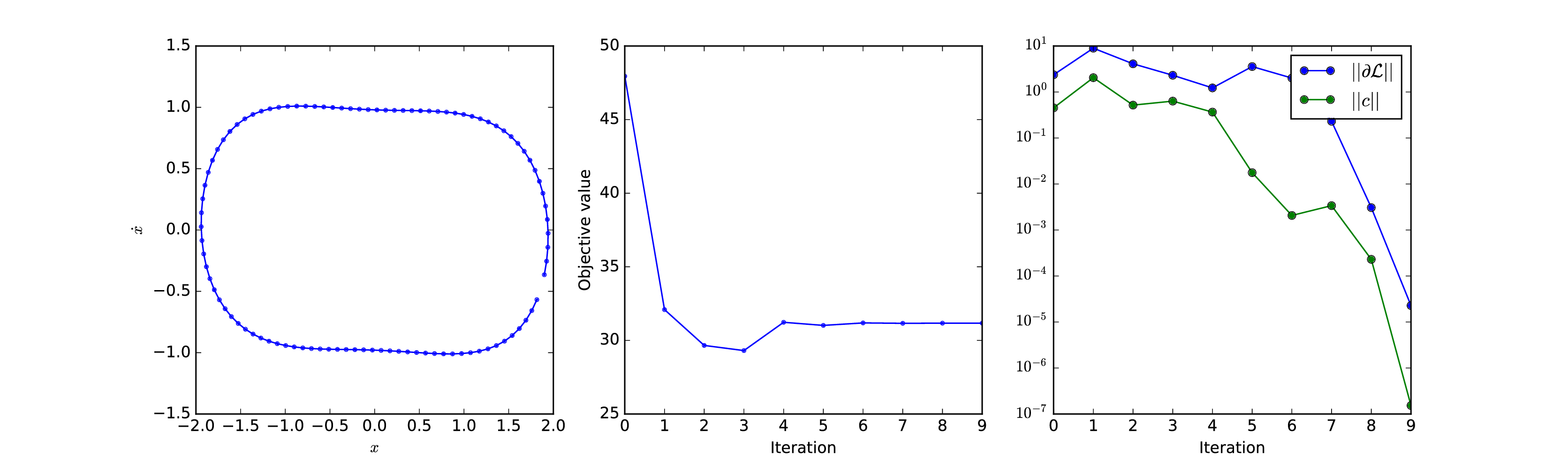}
\caption{\textsc{left}: An optimal limit cycle for the system \(\ddot{x} = -(x^3 + \dot{x}^3)/6 + u\). \textsc{middle}: Convergence of the objective function for the limit-cycle problem \eqnref{eqn:sdd}. \textsc{right}: Convergence in norm, of the Lagrangian gradient, and constraint deviation.}
\end{figure*}

\begin{equation}
  \label{eqn:mpcore}
\begin{tikzpicture}
  \matrix [matrix of math nodes,left delimiter={[},right delimiter={]}] (m)
  {
    \tilde{Q}_{\xi \xi} & \tilde{Q}_{\iota \xi}^{\mathtt{T}} & \tilde{G}_{:r, \xi}^{\mathtt{T}} & \tilde{G}_{r:, \iota}^{\mathtt{T}}\\
    \tilde{Q}_{\iota \xi} & \tilde{Q}_{\iota \iota} & \tilde{G}_{:r, \iota}^{\mathtt{T}}& \tilde{G}_{r:, \iota}^{\mathtt{T}}\\
    \tilde{G}_{:r, \xi} & \tilde{G}_{:r, \iota} & 0 & 0\\
    \tilde{G}_{r:, \xi} & \tilde{G}_{r:, \iota} & 0 & 0\\
  };
  \draw[color=red] (m-2-2.north west) -- (m-2-3.north east) -- +(0, -1.3) -- +(-1.75, -1.3) -- (m-2-2.north west);  
  \matrix [matrix of math nodes,left delimiter={[},right delimiter={]}] (v) at (3, 0)
  {
    \pr{\xi}\\
    \pr{\iota}\\
    \lambda_{:r}\\
    \lambda_{r:}\\
  };
  \node at (4, 0) {$=$};
  \matrix [matrix of math nodes,left delimiter={[},right delimiter={]}] (v) at (5, 0)
  {
    \tilde{b}_{\xi}\\
    \tilde{b}_{\iota}\\
    \tilde{h}_{:r}\\
    \tilde{h}_{r:}\\
  };
\end{tikzpicture}
\end{equation}

% \begin{equation}
%   \label{eqn:mpcore}
%   \left[\begin{array}{c c c c}
%         \tilde{Q}_{\xi \xi} & \tilde{Q}_{\iota \xi}^{\mathtt{T}} & \tilde{G}_{:r, \xi}^{\mathtt{T}} & \tilde{G}_{r:, \iota}^{\mathtt{T}}\\
%         \tilde{Q}_{\iota \xi} & \tilde{Q}_{\iota \iota} & \tilde{G}_{:r, \iota}^{\mathtt{T}}& \tilde{G}_{r:, \iota}^{\mathtt{T}}\\
%         \tilde{G}_{:r, \xi} & \tilde{G}_{:r, \iota} & 0 & 0\\
%         \tilde{G}_{r:, \xi} & \tilde{G}_{r:, \iota} & 0 & 0\\
%       \end{array}\right] \left[\begin{array}{c}
%                                  x_l\\
%                                  x_p\\
%                                  \lambda_{:r}\\
%                                  \lambda_{r:}\\
%                                \end{array}\right]  = \left[\begin{array}{c}
%                                                              \tilde{b}_{\xi}\\
%                                                              \tilde{b}_{\iota}\\
%                                                              \tilde{h}_{:r}\\
%                                                              \tilde{h}_{r:}\\
%                                                            \end{array}\right]
% \end{equation}

Gaussian Belief-Propagation is essentially a restatement of LU decomposition \cite{davis2006direct}. Gaussian-BP consists of messages of the form \cite{yedidia2000generalized} \cite{wainwright2008graphical},
\begin{equation}
  \begin{aligned}
    \mu_{i \rightarrow j} &:= [J_{i \rightarrow j}, h_{i \rightarrow j}] = [J_{ii}, h_i] - \sum_{k \in \delta(i) \backslash j} J_{ik} J_{k \rightarrow i}^{-1} [J_{ki}, h_{k \rightarrow i}],\\
    \mu_i & = J_{i \rightarrow j}^{-1} (h_{i \rightarrow j} - J_{i j} \mu_j),
  \end{aligned}
\end{equation}
where \(J \mu = h\) is the equation that is to be solved. These can be replaced by appropriate square-root forms to obtain instead, an LDL decomposition.

\begin{theo}
\label{theo:twlem}
The linear equation (\ref{eqn:spfamily}) can be solved in time \(\tilde{O}(\tw(\mathcal{H})^3)\), given the minimal tree-decomposition \emph{via} Algorithm~\ref{alg:gqp}.
\end{theo}
\emph{Proof.} The correctness of the algorithm follows from Lemma~\ref{lem:esep}. The bound holds trivially if, \(\rank{\tilde{G}} \leq \rank{\tilde{Q}_{\iota, \iota}}\), at every step of the algorithm. Otherwise, by realizing that \(\tilde{G}_{l \rightarrow p}\), can't have rank more than \(|\chi(p)|\), the proof follows.
\qed
%%Figure

%Figure
% Each individual message can also be derived from QPs and duality. Bertsekas in \cite{bertsekas1996thevenin} gives a \emph{feedback} interpretation of the algorithm by making use of the Legendre transform to solve the QP associated with (\ref{eqn:spfamily}).

It follows from Theorem~\ref{theo:twlem}, that the KKT system in Algorithm~\ref{alg:gnewt} can be solved in time \(\tilde{O}(\tw(\hat{\mathcal{G}})^3)\), where \(\hat{\mathcal{G}}\) is the moralization of the computational graph \(\mathcal{G}\).

The above proof also ensures that the equivalent sparse LU/LDL decomposition \cite{davis2006direct}, with the same pivot order, also has the same run-time. Since decompositions of indefinite systems are subject to instability, use of specialized solvers is generally preferable.

\section{Numerical Experiments}
In this section, we present preliminary numerical results with an implementation of Algorithm~\ref{alg:gnewt}, using the MA57 solver \cite{duff2004ma57}. For ensuring convergence in constrained problems, an augmented Lagragian merit function was used \cite{gill1986some}. The implementation was tested on the following non-standard control problems.

\begin{figure}[h]
\label{fig:sdgraph}
\includegraphics[height=5cm]{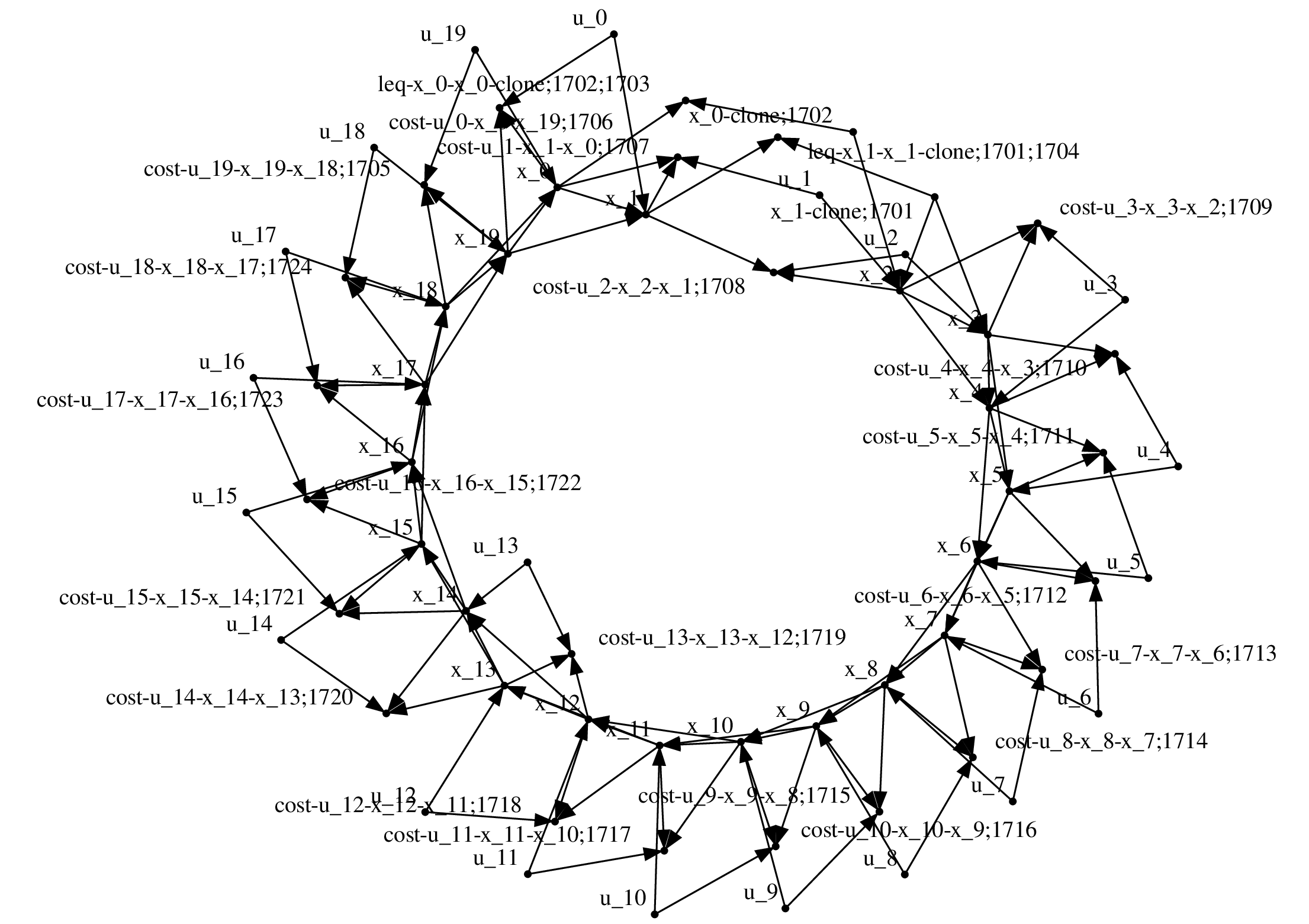}
\caption{The computational graph for the limit cycle problem, with \(N=20\)}
\end{figure}

\textsc{Spring-damper limit cycle}: Consider the following spring-damper limit-cycle problem \cite{tassa2011optimal} (Figure~\ref{fig:sdgraph}),
\begin{equation}
\begin{aligned}
&\min_{x_0, u[0, T]} \int \ell(x, \dot{x}, u) \diff{t},\\
&\ddot{x} = -(x^3 + \dot{x}^3)/6 + u,\\
&x(0) = x(T) = x_0, \dot{x}(0) = \dot{x}(T) = \dot{x}_0,
\end{aligned}
\end{equation}
where,
\[\ell(x, \dot{x}, u) = (1 - \e^{(\dot{x}_i - 2)^2} - \e^{- (\dot{x}_i + 2)^2}) + {1 \over 2} \norm{u_i}^2_2.\]

Discretising the derivatives by finite differences, \(\dot{x} \approx \Delta x_i / \Delta t = (x_i - x_{i - 1})/\Delta t\), this can be written as the following structured optimization problem,
\begin{equation}
\label{eqn:sdd}
\begin{aligned}
& \quad\quad \min_{x_0, \{u_i\}_{0}^m} \sum_{1}^N \ell(x_i, \Delta x_i/\Delta t, u_i),\\
&x_{i + 1} \leftarrow x_{i} + \Delta x_{i} + (\Delta t)^2 \left[ -(x^3 + (\Delta x_{i} / \Delta t)^3)/6 + u_i\right].\\
&x_{0} = x_{N - 2}, x_{1} = x_{N - 1}.
 % \quad x(0) = x(T) = x_0,
\end{aligned}
\end{equation}

For \(N = 100, \Delta t = 0.1\), with random initializations, the problem showed robust convergence; often taking no more than ten SQP iterations. The optimal limit cycle, and the convergence curves for one run of the algorithm are shown in (Figure~\ref{fig:sd}).

\begin{figure}[h!]
\label{fig:ab}
\centering
\includegraphics[height=10cm]{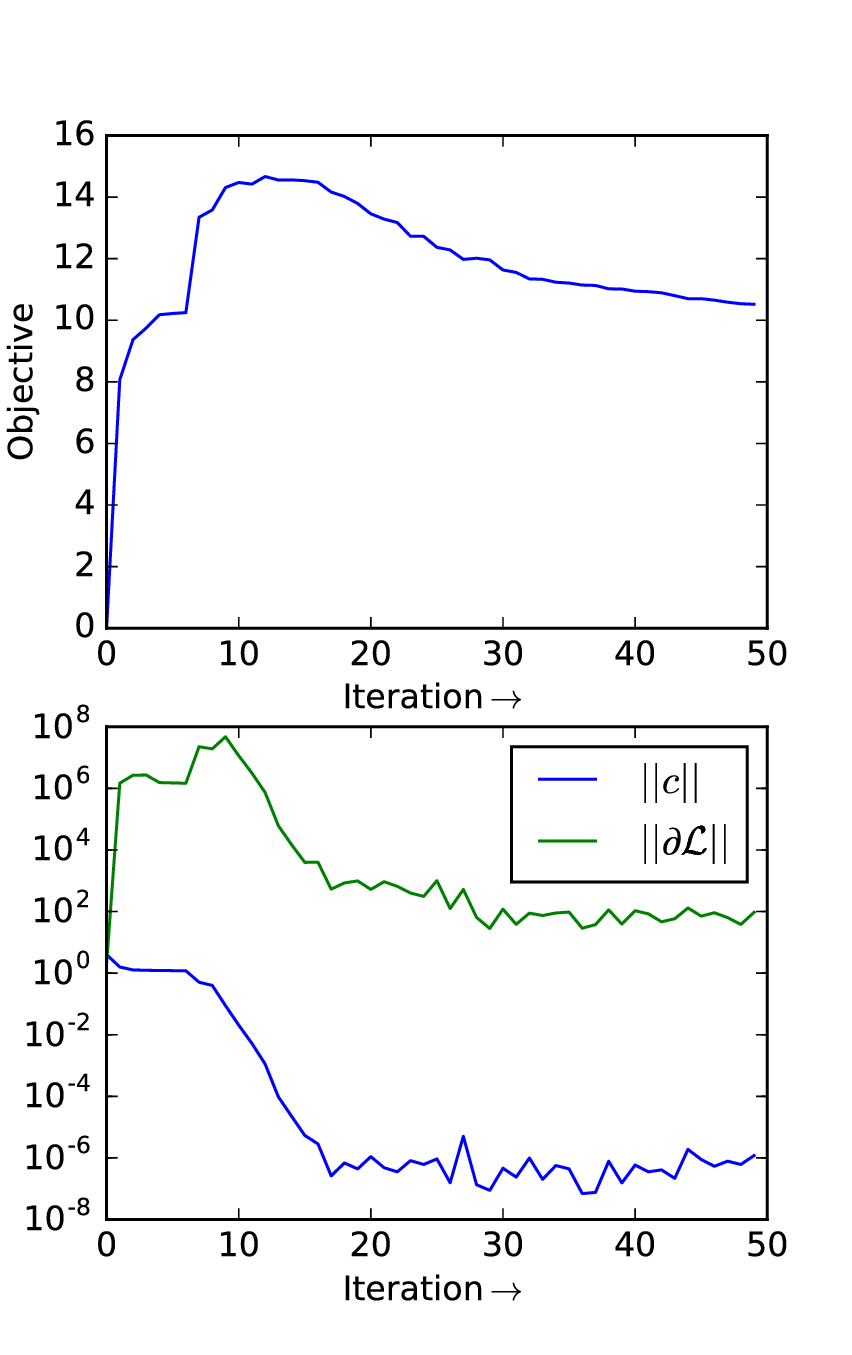}
\caption{\textsc{top}: Convergence of the objective function for the Acrobot problem. \textsc{bottom}: Convergence in norm, of the Lagrangian gradient, and constraint deviation.}
\end{figure}

\textsc{Acrobot}:
We also considered the Acrobot problem, using a discretized second-order dynamics as before, and \(\ell_2\) and total-variation penaltes over the control sequence; the end-condition being enforced using a constraint over the final cartesian position. For \(N = 100, \Delta t = 0.04d0\), with random initializations, the optimizer converged often in about 20 iterations (Figure~\ref{fig:ab}).

% \begin{equation}
%   \begin{aligned}
%     &\int \ell(x, \dot{x}, u) \diff{t}
%     &\sum_{} \ell(x_i, x_{i - 1}, u_i) \\
%     & \ell(x, \dot{x}, u) = (1 - \e^{(\dot{x}_i - 2)^2} - \e^{- (\dot{x}_i + 2)^2}) + {1 \over 2} \norm{u_i}^2_2.
%   \end{aligned}
% \end{equation}

% \begin{equation}
% \begin{aligned}
% &\sum_{} \ell(x_i, x_{i - 1}, u_i) \\
% & \ell(x, \dot{x}, u) = (1 - \e^{(\dot{x}_i - 2)^2} - \e^{- (\dot{x}_i + 2)^2}) + {1 \over 2} \norm{u_i}^2_2.
% \end{aligned}
% \end{equation}

\section{Discussion}
We have shown that the Newton step can be computed in time \(\tilde{O}(\tw^3)\), where '\(\tw\)' is the tree-width of the computational graph. We have also derived extensions to constrained problems, and provided numerical examples. The technique presented herein, also generalizes many specialized algorithms in control.

In certain control problems, the solution to the KKT system, itself can be written in \emph{feedback form}. Given a \(LU\) decomposition of the KKT system, one can replace the backsubstitution phase by \(U\), with a function evaluation that uses \(L\) as a \emph{control feedback} \cite{jacobson1970differential}. It is unclear if such techniques can be generalized, and whether they can be made independent of the pivot-order used for solving the system.

A competing method for exploiting the structure of objectives such as \eqnref{eqn:obj}, is by the use Hessian vector product AD routines in conjugation with CG-like methods. Computing the Hessian vector product takes time \(\tilde{O}(\omega(\hat{\mathcal{G}})^2)\), where \(\hat{\mathcal{G}}\) is a moralization of the computational graph \cite{dixon2009automatic}. By contrast, if the computational graph were chordal, then computing the Newton-step \emph{via} Algorithm~\ref{alg:gnewt} is only \(\tilde{O}(\omega(\hat{\mathcal{G}})^3)\). The latter is more economical when the cliques of a graph are small in comparison to the order of the graph. The ill-conditioned nature of structured objectives may also lead to bad convergence properties for CG algorithms.

For problems whose tree-widths are large, the iterative method is obviously more viable. However, following the rapid advances in approximate inference in the past two decades \cite{wainwright2008graphical}, we hope that the explicit algebraic connection to graphical models made in this paper, can be exploited in coming up with less-agnostic iterative methods.


%@Comment author={Srinivasan, Akshay and Todorov, Emanuel},

\bibliographystyle{IEEEtran}
\bibliography{references}

\end{document}